\documentclass[notitlepage]{jams-l} 
\usepackage{amssymb,amsmath,latexsym,epsfig,euscript,floatfig,wrapfig} 

\newlength{\jmr}
\newlength{\hwl}
\newlength{\khov}
\newtheorem{smi}{Smirnov's Theorem} 

\newtheorem{nota}{Notation} 

\newtheorem{lip}{Lipshitz's Theorem}

\newtheorem{thm}{Theorem}
\newtheorem{lemma}{Lemma}

\newtheorem{dfn}{Definition}
\newtheorem{cor}{Corollary}
\newtheorem{rem}{Remark}	
\newtheorem{ex}{Example}

\newcommand{\eps}{\varepsilon}

\newcommand{\supp}{\mathrm{Supp}}

\newcommand{\conv}{\mathrm{Conv}}

\newcommand{\thth}{{\underline{\mathrm{th}}}}

\newcommand{\nd}{{\underline{\mathrm{nd}}}}

\newcommand{\F}{\mathbb{F}}
\newcommand{\Q}{\mathbb{Q}}
\newcommand{\R}{\mathbb{R}}
\newcommand{\C}{\mathbb{C}}
\newcommand{\N}{\mathbb{N}}
\newcommand{\Z}{\mathbb{Z}}
 
\newcommand{\cp}{\mathfrak{p}}

\newcommand{\ord}{\mathrm{ord}}

\newcommand{\Zn}{\Z^n}

\newcommand{\Rn}{\R^n}

\renewcommand{\qed}{$\blacksquare$}
\newcommand{\cL}{{\mathcal{L}}}
\newcommand{\cM}{{\mathcal{M}}}

\newcommand{\bO}{\mathbf{O}}

\begin{document}

\title[Arithmetic Fewnomial Systems]{Finiteness for 
Arithmetic Fewnomial Systems} 


\author{J.\ Maurice Rojas}\thanks{
This research was partially supported by Hong Kong UGC Grant \#9040469-730 
and a grant from the Texas A\&M Faculty of Sciences.} 

\address{Department of Mathematics, Texas A\&M University, College Station, 
Texas 77843-3368, USA. }  
\email{rojas@math.tamu.edu } 
\urladdr{http://www.math.tamu.edu/\~{}rojas } 

\dedicatory{To my god-daughter, Monica Althoff. }  

\date{\today} 

\mbox{}\\
\vspace{-.9in} 
\maketitle

\mbox{}\\
\vspace{-1.5cm} 
\begin{abstract} 
Suppose $\cL$ is any finite algebraic extension of either the ordinary 
rational numbers or the $p$-adic rational numbers. Also let $g_1,\ldots,g_k$ 
be polynomials in $n$ variables, with coefficients in $\cL$, such that the 
total number of monomial terms appearing in at least one $g_i$ is exactly $m$. 
We prove that the maximum number of isolated roots of 
$G\!:=\!(g_1,\ldots,g_k)$ in $\cL^n$ 
is finite and depends solely on $(m,n,\cL)$, i.e., is independent of the 
degrees of the $g_i$. We thus obtain an arithmetic analogue of Khovanski's 
Theorem on Fewnomials, extending earlier work of Denef, Van den Dries, 
Lipshitz, and Lenstra. 
\end{abstract} 

\section{Introduction} 
A consequence of Descartes' Rule (a classic result dating back to 1641) 
is that any real univariate polynomial with exactly $m\!\geq\!1$ monomial 
terms has at most $2m-1$ real roots. In this paper, we extend this result 
in two directions at once: we bound the number of isolated roots of 
polynomial {\bf systems} over {\bf $\pmb{\cp}$-adic} fields, independent of 
the degrees of the underlying polynomials. As 
a consequence, we also obtain analogous results over any number 
field. The resulting bounds are, unfortunately, non-explicit. 
Nevertheless, the existence of such bounds on the number of roots was 
previously unknown in the multivariate non-Archimedean and number field cases. 
So let us now detail our results and some important precursors. 

Extending the sharp bound of Descartes' Rule to 
polynomial systems has already proven difficult in the special case 
of the real numbers $\R$ (perhaps the most familiar metrically complete field): 
the best general result in this direction, Khovanski's Theorem on Fewnomials  
\cite{kho,few}, gives an explicit bound on the number of non-degenerate 
real roots\footnote{One can then extend this to counting isolated roots and 
connected components via some algebraic tricks. See, e.g., \cite{real,tri}.} 
that is independent of the underlying polynomial degrees,  
but exponential in the number of monomial terms. 
Whether the number of non-degenerate roots is in fact {\bf polynomial} 
in the number of monomial terms for every {\bf fixed} number of 
equations and variables is an intriguing open question, and the 
answer is still unknown even in the case of two polynomials in two unknowns. 
(The polynomial system $(x^2_1-x_1,\ldots,x^2_n-x_n)$ easily shows us 
that fixing $n$ is necessary if we would like polynomiality in the number of 
monomial terms.) 

As for the $p$-adic univariate case, Hendrik W.\ Lenstra, Jr.\ has shown that 
for any degree $d$ algebraic extension $L$ of $\Q_p$, the maximum number of 
roots in $L$ of a $g\!\in\!L[x_1]$ with exactly $m\!\geq\!2$ 
monomial terms is no more than $1+1.582 \cdot (p^{f_L}-1)(m-1)^2(1+
\frac{e_L\log(e_L(m-1)/\log p)}{\log p})$ \cite[prop.\ 7.2]{lenstra}, 
where $e_L$ and $f_L$ respectively denote the {\bf ramification degree} and 
{\bf residue field degree} of $L$. (Recall that $e_L$ and $f_L$ are integers 
satisfying $e_Lf_L\!=\!d$; see section \ref{sec:local} for their definitions.) 
In particular, this bound is independent of the degree of $g$. As a 
consequence, Lenstra also derived a bound of $1+4.566\cdot 
(m-1)^2(d+10)2^d(\log(d(m-1))+0.367)$ \cite[prop.\ 8.1]
{lenstra} for the analogous situation where one replaces $L$ by any degree 
$d$ algebraic extension $K$ of $\Q$.\footnote{
In both cases, a bound of $1$ is 
trivial to derive when $m\!\leq\!1$, and the bounds we quoted have the added 
benefit of counting the multiplicities of the non-zero roots. Furthermore, the 
result over $K$ actually is a stronger statement which in fact bounds the 
number of roots of bounded degree over $K$. } 

By recent work of the author \cite{tri} this polynomiality  
of the number of roots in the number of monomials can be extended 
to certain systems of $n$ polynomials in $n$ unknowns, provided 
we fix $n$ and restrict to real algebraic number fields. 
(The example from the paragraph before last tells us that 
fixing $n$ is necessary in the $\cp$-adic and number field cases as well.) 
However, at the expense of less explicit bounds, one can 
extend Lenstra's results much farther. 

\begin{nota} 
Let $\mathbf{\cL}$ be a field and $\cL^*\!:=\!\cL\setminus\{0\}$. 
If $\pmb{G}\!:=\!(g_1,\ldots,g_k)$ where, for all 
$i$, \scalebox{.85}[1]{$\pmb{g_i}\!\in\!\cL[x_1,\ldots,x_n]\setminus\{0\}$,} 
and the number of monomial terms appearing in at least one $g_i$ is exactly 
$m$, then we call $G$ a {\bf ($\pmb{k\times n}$)  
$\pmb{m}$-sparse polynomial system (over $\pmb{\cL}$)}. 
Also, we say a root $\zeta$ of $G$ is {\bf isolated} 
(resp.\ {\bf non-degenerate}) iff $\zeta$ is an irreducible component 
of the zero set of $G$ over the algebraic closure of $\cL$ 
(resp.\ $k\!=\!n$ and the Jacobian of $G$, evaluated at $\zeta$, is 
invertible). 
\end{nota} 

\begin{thm}
\label{thm:local}
For any (rational) prime $p$ and positive integer $d$, let $L$ be any degree 
$d$ algebraic extension of $\Q_p$. Also let $G$ be any 
$k\!\times\!n$ $m$-sparse polynomial system over $L$. Then there is an 
absolute constant $\gamma(n,m)$ such that the number of isolated roots of $G$ 
in $L^n$ is no more \mbox{than $p^{dn}(1-\frac{1}{p^{f_L}})^n\gamma(n,m)$.} 
\end{thm} 
\begin{cor} 
\label{cor:global}
Let $K$ be any degree $d$ algebraic extension of $\Q$ and let $G$ be 
any $k\!\times\!n$ $m$-sparse polynomial system over $K$. Then 
the number of isolated roots of $G$ in $K^n$ 
is no more \mbox{than $2^{dn}(1-\frac{1}{2^d})^n\gamma(n,m)$.} 
\end{cor} 
\noindent 
Theorem \ref{thm:local} generalizes an analogy over $\Z_p$, 
initiated by Jan Denef and Lou Van den Dries in \cite{vandenef}, of 
Khovanski's Theorem on Fewnomials. Corollary \ref{cor:global} establishes 
a higher-dimensional analogue of Lenstra's aforementioned result for 
univariate sparse polynomials over number fields. We can also extend 
our finiteness results even further to count isolated roots of bounded degree 
over $L$ or $K$ (see corollary \ref{cor:bigger} of section \ref{sec:extend}). 

We prove theorem \ref{thm:local} and corollary \ref{cor:global}
in sections \ref{sec:local} and \ref{sec:global} respectively.  
The proofs, while short, involve deep non-effective results of Jan 
Denef and Lou Van den Dries \cite{vandenef} and Leonard Lipshitz 
\cite{lipshitz} on $p$-adic sub-analytic functions,  
as well as an elegant extension of the 
classical $p$-adic Newton polygon by A.\ L.\ Smirnov \cite{smirnov}. In 
particular, aside from the case $n\!=\!1$ (cf.\ remark \ref{rem:heart} of 
the next section), there appear to be no explicit bounds on the function 
$\gamma(n,m)$ yet. So a more direct and effective approach would be of the 
utmost interest. 

\subsection{$p$-adic Analysis and $p$-adic Newton Polytopes} 
\label{sec:padic} 
\mbox{}\\

We first state the following combined paraphrase of two results of Lipshitz: 
\begin{lip} (See \cite[thm.\ 2]{lipshitz}.) 
For any (rational) prime $p$, let $\C_p$ denote the completion (with respect 
to any $p$-adic metric) of the algebraic closure of $\Q_p$. 
Also let $G$ be any $n\!\times\!n$ $m$-sparse polynomial system over $\C_p$. 
Then there is an absolute constant $\beta'(n,m)$ (independent of $p$) such 
that $G$ has no more than $\beta'(n,m)$ isolated roots 
$x\!:=\!(x_1,\ldots,x_n)\!\in\!\C^n_p$ 
satisfying $|x_i-1|_p\!\leq\!\frac{1}{p}$ for all $i$, where $|\cdot|_p$ 
denotes the unique $p$-adic norm on $\C_p$ with $|p|_p\!=\!\frac{1}{p}$. 
\qed
\end{lip} 
\begin{rem} 
In the above $p$-adic context, we also have the following equivalent 
definition of isolation for roots: a root 
$x$ of $G$ is isolated iff for some $\eps\!>\!0$, 
we have $\max_i|x_i-x'_i|_p\!>\!\eps$ for every other root 
$(x'_1,\ldots,x'_n)$ of $G$. So in essence, an isolated root of $G$ in 
$\C^n_p$ can be contained within a small $p$-adic ``brick,'' away from all 
other roots of $G$. Lipshitz's original statement in fact dealt with roots with 
{\bf algebraic integer} coordinates in $\C^n_p$, but the statement above is 
equivalent since the ultrametric inequality implies 
$|x|_p\!\leq\!\max\{|x-1|_p,|1|_p\}\!=\!1$.  
\end{rem} 
\begin{rem} 
\label{rem:heart} 
Lenstra has derived an explicit 
upper bound on the number of roots $x_1$ of $g_1$ in $\C_p$ with 
$|x_1-1|_p\!\leq\!\frac{1}{p^{r}}$ for any given 
$r\!>\!0$ \cite[prop.\ 7.1]{lenstra}.\footnote{ The original statement was in 
terms of $\ord_p$, counted 
multiplicities, and in fact gave a decreasing function of $p$. } 
Taking $r\!=\!1$ one then obtains 
$\beta'(1,m)\!\leq\!1.582\cdot(m-1)(1+1.443\cdot(\log(m-1)+0.367))$ for 
$m\!\geq\!2$. (Note that $\beta'(n,0)\!=\!\beta'(n,1)\!=\!0$ for all $n$.) 
Whether Lenstra's explicit bound on the number of roots in $\C_p$ 
``$p$-adically close to the identity'' extends to $n\times n$ sparse 
polynomial systems is an open problem, even in the case $n\!=\!2$. 
Nevertheless, the proofs of theorem \ref{thm:local} and corollary 
\ref{cor:global} are structured so that explicit bounds on $\gamma(n,m)$ can 
be easily derived should such a result become available. 
\end{rem} 

\noindent 
Lipshitz's Theorem is based partially on an earlier result of 
Denef and Van den Dries \cite[pg.\ 105]{vandenef} over the subring $\Z_p$ 
but also injects model-theoretic techniques (see \cite{lipshitz} for further 
details).  

The key to proving theorem \ref{thm:local} is to further limit 
the number of roots defined over a subfield of $\C_p$ by seeing which 
possible vectors of valuations can occur. In particular, we will 
use the following extension of the classical univariate 
$p$-adic Newton polygon (see, e.g., \cite[ch.\ IV, sec.\ 3]{koblitz} for the 
latter construction). To clarify the statement, let us make the following 
definitions: 
\begin{dfn} 
\label{dfn:poly}
For any $a\!=\!(a_1,\ldots,a_n)\!\in\!\Zn$, let 
$x^a\!:=\!x^{a_1}_1\cdots x^{a_n}_n$. 
Writing any $g\!\in\!\cL[x_1,\ldots,x_n]$ in the form $\sum_{a\in \Zn} 
c_ax^a$, we call $\pmb{\supp(g)}\!:=\!\{a \; | \;
c_a\!\neq\!0\}$ the {\bf support} of $g$. 
Then, for any $n\times n$ polynomial system $G$ over $\C_p$, its {\bf 
$\pmb{n}$-tuple of $\pmb{p}$-adic Newton polytopes}, 
$\pmb{\Delta_p(G)}\!=\!(\Delta_p(g_1),\ldots,\Delta_p(g_n))$, is 
defined as follows: $\pmb{\Delta_p(g_i)}\!:=\!\conv(\{(a,\ord_p c_a) \; | \; 
a\!\in\!\supp(g_i) \})\!\subset\!\R^{n+1}$, where $\conv(S)$ denotes the 
convex hull of\footnote{i.e., smallest convex set containing...} 
a set $S\!\subseteq\!\R^{n+1}$ and $\ord_p : \C_p \longrightarrow 
\Q\cup\{+\infty\}$ is the usual discrete valuation\footnote{So, for example, 
$\ord_p 0\!=\!+\infty$ and $\ord_p(p^kr)\!=\!k$ whenever $r$ is a unit in 
$\Z_p$ and $k\!\in\!\Q$.} of $\C_p$. Finally, for any $w\!\in\!\Rn$ and any 
compact subset $B\!\subset\!\Rn$, let the {\bf face of $\pmb{B}$ with inner 
normal $\pmb{w}$}, $\pmb{B^w}$, be the set of points $x\!\in\!B$ which 
minimize the inner product $w\cdot x$. 
\end{dfn} 
\begin{ex}
Consider the $2\!\times\!2$ $6$-sparse polynomial system 
\[F\!:=\!(50x^{18}_1-3125x^9_2-162x_2,49x^{18}_2-35x^9_1-109375x_1)\] 
over $\Q_5$. Then the corresponding pair of $5$-adic Newton 
polytopes is  
\[\Delta_5(F)\!=\!(\conv(\{(18,0,2),(0,9,5),(0,1,0)\}),
\conv(\{(0,18,0),(9,0,1),(1,0,6)\})).\] 
Note that each polytope is in fact a triangle embedded in $\R^3$. 
\end{ex}  

\begin{smi} 
\cite[thm.\ 3.4]{smirnov}
\label{thm:smirnov}
Let $\ord_p x$ be the vector $(\ord_p x_1,\ldots,\ord_p x_n)$ and 
let $v\!:=\!(v_1,\ldots,v_n)\!\in\!\Rn$. Then 
for any $n\!\times\!n$ polynomial system $G$ over $\C_p$, 
the number of isolated roots $x\!:=\!(x_1,\ldots,x_n)$ of $G$ in $(\C^*_p)^n$ 
satisfying $\ord_p x\!=\!v$ (counting multiplicities) 
is no more than $\cM(\pi(\Delta^{v'}_p(g_1)),\ldots,\pi(\Delta^{v'}_p(g_n)))$, 
where $v'\!:=\!(v_1,\ldots,v_n,1)$, $\pi : \R^{n+1}\longrightarrow \Rn$ is the 
natural projection forgetting the $v_{n+1}$ coordinate, $\cM$ denotes 
mixed volume \cite{buza} (normalized so that 
\scalebox{1}[1]{$\cM(\conv(\{\bO,e_1,\ldots,e_n\}),
\ldots,\conv(\{\bO,e_1,\ldots,e_n\}))\!=\!1$}), and 
$e_i$ is the $i^\thth$ standard basis vector of $\Rn$. \qed 
\end{smi} 
\begin{rem} 
\label{rem:convex} 
Note that the number of roots of $G$ in $(\C^*_p)^n$ with given 
valuation vector 
thus depends strongly on the individual exponents of $G$ --- not just on 
the number of monomial terms. However, the number of possible distinct 
valuation vectors occuring for any single $G$ can be combinatorially bounded 
from above as a function depending solely on $n$ and the number of monomial 
terms (cf.\ section \ref{sec:local}). In particular, it is only the 
{\bf lower}\footnote{Those with positive 
$x_{n+1}$ coordinate for their inner normals...} faces of the Newton polytopes 
that matter. 
\end{rem} 

\noindent 
Note that Smirnov's Theorem provides a non-Archimedean extension of 
Bernstein's famous mixed volume bound over $(\C^*)^n$ \cite{bkk,jpaa}.  
We also point out that aside from a result of Kamel A.\ Atan and 
J.\ H.\ Loxton in the $2\times 2$ case \cite{loxton}, Smirnov's result appears 
to be the first higher-dimensional version of the classical univariate  
$p$-adic Newton polygon. 

While we will leave the algorithmic issues of $p$-adic Newton polytopes 
for another paper, let us at least observe one salient fact before pointing 
out references to the computational literature: In searching for 
$v'\!\in\!\R^{n+1}$ giving a positive number of roots with valuation vector 
$v$, it suffices to restrict one's search to the inner normals of the lower 
$n$-dimensional faces of the {\bf Minkowski sum}\footnote{The Minkowski sum 
of any finite collection of subsets $A_1,\ldots,A_k\subseteq\!\Rn$ is simply 
the set $\{a_1+\cdots+a_k \; | \; a_i\in A_i$ for all $i\}$. } 
$\Delta_p(g_1)+\cdots+\Delta_p(g_n)$. The last fact follows from basic convex 
geometry (see, e.g., \cite{buza}). While there currently seems to 
be no direct software implementation of $p$-adic Newton polytopes, the 
underlying algorithms have already been implemented in the related context 
of mixed volume computation, and a detailed description including 
complexity bounds can be found in \cite{emiphd,sparse}. 
\begin{ex}
A convenient way to visualize how many roots with valuation vector $v$ 
appear for a given $2\times 2$ polynomial system is to draw the projected 
lower faces of the underlying Minkowski sum. For instance, in our last 
example, we obtain the following representation:\footnote{In the 
illustration, the {\bf ordinary} Newton polygons we refer to are 
a construction similar to the $p$-adic Newton polygon, embedded in $\R^2$ 
instead of $\R^3$, where one essentially uses the {\bf trivial} valuation 
($|a|\!=\!1$ for any $a\!\in\!\C^*_p$) instead of the $p$-adic 
valuation. }\\ 
\mbox{}\hspace{1in}\epsfig{file=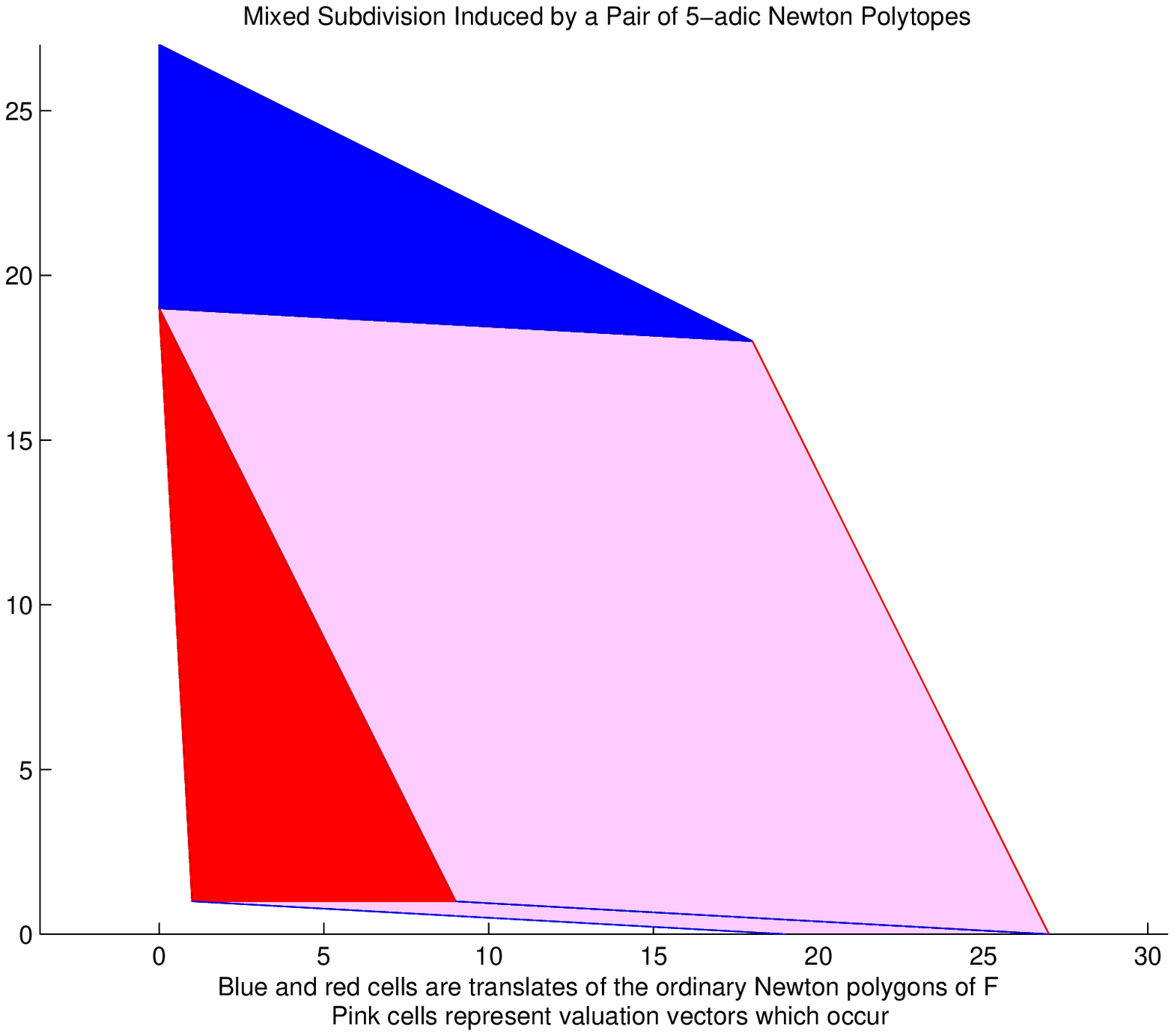,height=2.9in} 

It is then easily checked\footnote{In this case, a brute-force search 
among the cross-products of the pairs of triangle edges suffices to generate 
our normal vectors: $2$ out of the resulting $9$ possibilities are true 
inner normals of lower $2$-dimensional faces of the underlying Minkowski sum. 
Each resulting pair $(\pi(\Delta^{v'}_5(g_1)),\pi(\Delta^{v'}_5(g_2)))$ 
turns out to be a pair of line segments. The mixed area of any such pair 
(whose Minkowski sum is a lighter-colored cell in our illustration) is simply 
the absolute value of the determinant of the underlying vectors. }   
that among the $323$ roots $(x_1,x_2)$ of $G$ in $(\C^*_5)^2$, there are 
exactly $8$ with $(\ord_5 x_1,\ord_5 x_2)\!=\!(\frac{5}{8},\frac{53}{4})$ 
and $315$ with $(\ord_5 x_1,\ord_5 x_2)\!=\!(-\frac{1}{9},0)$. 
In particular, we see that there are roots lying in extensions of $\Q_5$ 
of degree at least $9$. 
\end{ex} 

\section{Roots of Bounded Degree Over $\cp$-adic Fields 
and Number Fields } 
\label{sec:extend} 
Here we use an observation on $p$-adic algebraic extensions to 
prove the following combined strengthening of theorem \ref{thm:local} 
and corollary \ref{cor:global}. First, let us say that 
a vector $x\!:=\!(x_1,\ldots,x_n)$ defined over the algebraic closure of a 
field $\cL$ {\bf is of degree $\pmb{\delta}$ over $\pmb{\cL}$} iff each 
$x_i$ is of degree $\delta$ over $\cL$. 
\begin{cor}
\label{cor:bigger}  
Suppose $G$ is a $k\!\times\!n$ $m$-sparse polynomial system over $\cL$, 
where $\cL$ is a degree $d$ algebraic extension of either $\Q$ or $\Q_p$. 
Then for any positive integer $\delta$, there is an absolute constant 
$\gamma'(d,\delta,n,m)$ (resp.\ 
$\gamma'_p(d,\delta,n,m)$) such that $G$ has no more than 
$\gamma'(d,\delta,n,m)$ (resp.\ 
$\gamma'_p(d,\delta,n,m)$) isolated roots in 
$\C^n$ (resp.\ $\C^n_p$) 
of degree $\leq\!\delta$ over $\cL$, 
according as $\cL$ is an algebraic extension of $\Q$ or $\Q_p$.   
\end{cor} 

\noindent 
{\bf Proof:} Focusing first on the case where $\cL$ is an algebraic extension 
of $\Q_p$, note that there are only finitely many algebraic extensions 
of degree $\leq\!d\delta$ of $\Q_p$ \cite[ch.\ II, prop.\ 14]{lang}. 
Letting $L$ be the compositum of all these fields, note that 
$L$ is then also a finite algebraic extension of $\Q_p$.  
More to the point, any root of $G$ in $\C^n_p$ of degree 
$\leq\!\delta$ over $\cL$ must then also lie in $L^n$. 
So the $p$-adic case of our corollary follows 
immediately from theorem \ref{thm:local}. 

To prove the case where $\cL$ is an algebraic extension of 
$\Q$, note that such an $\cL$ embeds naturally as a subfield of the $L$ we 
just defined for the $p$-adic case. Taking $p\!=\!2$ to fix ideas, 
we then see that the degree of $\cL\Q_2$ over $\Q_2$ is 
no more than $d$, and thus any $x_i$ of degree $\leq\!\delta$ over $\cL$ 
embeds in an extension of $\Q_2$ of degree $\leq\!d\delta$. Thus any such 
$x_i$ can be assumed to lie in $L$, and we again conclude by theorem 
\ref{thm:local}. \qed 

\begin{rem} 
Note that we immediately obtain from our proof above 
that \[\gamma'_p(d,\delta,n,m) \leq p^{nD_p}(1-\frac{1}{p^{f_p}})^n\gamma(n,m) 
\text{ \ and \ } 
\gamma'(d,\delta,n,m)
\leq 2^{nD_p}(1-\frac{1}{2^{f_2}})^n\gamma(n,m),\] 
where $D_p$ (resp.\ $f_p$) is the degree (resp.\ residue field degree) over 
$\Q_p$ of the compositum of all algebraic extensions of $\Q_p$ of degree 
$\leq\!d\delta$.  
\end{rem}  

\section{Proving our Main Local Result (Theorem \ref{thm:local})} 
\label{sec:local} 

The following lemma will allow us to reduce to the case 
$k\!=\!n$. 
\begin{lemma}
\label{lemma:gh}
Following the notation of theorem \ref{thm:local}, there is a matrix 
$[a_{ij}]\!\subset\!\Q^{n\times k}_p$ such that
the zero set of 
$G\!:=\!(a_{11}g_1+\cdots+a_{1k}g_k,\ldots,a_{n1}g_1+\cdots+a_{nk}g_k)$
in $\C^n_p$ is the union of the zero set of $G$ in $\C^n_p$ and a finite
(possibly empty) set of points. \qed  
\end{lemma}  
\noindent 
A stronger version of the above lemma appears 
in \cite[sec.\ 3.4.1]{giustiheintz}, but phrased over $\C$ instead. 
However, the proof there carries over to any algebraically closed field with 
no difficulty whatsoever. 

Returning to the proof of theorem \ref{thm:local}, we see that lemma 
\ref{lemma:gh} allows us to replace $G$ by a new 
$n\times n$ polynomial system (clearly still $m$-sparse) which 
has at least as many isolated roots as our original $G$. 
Abusing notation slightly, let $G$ denote this new $n\times n$ 
polynomial system. 

Applying Smirnov's Theorem to $G$, recall that 
$\cM(\pi(\Delta^{v'}_p(g_1)),\ldots,\pi(\Delta^{v'}_p(g_n)))\!>\!0 
\Longrightarrow \ v'$ is an inner normal of a lower $n$-dimensional face of 
the Minkowski sum $\Sigma_G\!:=\!\Delta_p(g_1)+\cdots+\Delta_p(g_n)$ (cf.\ 
section \ref{sec:padic}). 
It is then easily checked 
that $\Sigma_G$ has at most $mn$ vertices and thus, since 
any $n$-dimensional face consists of at least $n+1$ vertices, $\Sigma_G$ 
has at most $\begin{pmatrix}mn\\ n+1\end{pmatrix}$ 
$n$-dimensional faces. In particular, this implies that 
the number of distinct values for the vector 
$\ord_p x$, where $x\!\in\!(\C^*_p)^n$ is a root of $G$, is no more than 
$\begin{pmatrix}mn\\ n+1\end{pmatrix}$. So let us fix $v\!\in\!\Rn$ and 
see how many roots of $G$ in $(L^*)^n$ can have valuation vector $v$. 

Let $R_p\!:=\!\{a\!\in\!\C_p \; 
| |a|_p\!\geq\!1\}$ be the ring of algebraic integers in $\C_p$, $M_p$ the 
unique maximal ideal of $R_p$, $\F_L\!:=\!(R_p\cap L)/(M_p\cap L)$,  
and let $\pi$ be any generator of the principal ideal $M_p\cap L$ 
of $R_p\cap L$. Also let $e\!:=\!\max_{y\in L^*}\{ 
|\ord_p y|^{-1}\}$ and $f\!:=\!\log_p\#\F_L$.  
The last two quantities are respectively known as the {\bf ramification 
degree} and {\bf residue field degree} of $L$,  
and we can in fact pick $\pi$ so that $\pi^e\!=\!p$ as well 
\cite[ch.\ III]{koblitz}. Doing this, then 
fixing a set $A_L\!\subset\!R_p$ of representatives for $\F_L$ 
(i.e., a set of $p^f$ elements of $R_p\cap L$, exactly {\bf one} of 
which lies in $M_p$, whose image 
mod $M_p\cap L$ is $\F_L$), we can then write any 
$x_i\!\in\!L$ uniquely as $\sum^{+\infty}_{j=e\ord_p(x_i)}a^{(i)}_j\pi^j$ for 
some sequence of $a^{(i)}_j\!\in\!A_L$ \cite[cor., pg.\ 68, sec.\ 3, ch.\ III]
{koblitz}. 

Note in particular that 
$\frac{\sum^{+\infty}_{j=t}c_j\pi^j}{\pi^t(c_t+\cdots+c_{e+t-1}\pi^{e-1})}
\!\in\!R_p$, and in fact $\left|\frac{\sum^{+\infty}_{j=t}c_j\pi^j}
{\pi^t(c_t+\cdots+c_{e+t-1}\pi^{e-1})}-1\right|_p\!\leq\!\frac{1}{p}$,  
for any sequence of representatives $(c_t,c_{t+1},\ldots)\!\in\!A^{+\infty}_L$ 
with $c_t\!\in\!A_L\setminus M_p$. 

Now consider the polynomial system $H$ where  
\[H(z_1,\ldots,z_n)\!:=\!G(\pi^{ev_1}(c^{(1)}_{ev_1}+\cdots+c^{(1)}_{ev_1+e-1}
\pi^{e-1})z_1,\ldots,
\pi^{ev_n}(c^{(n)}_{ev_n}+\cdots+c^{(n)}_{ev_n+e-1}\pi^{e-1})z_n),\] 
for any fixed vectors $(c^{(1)}_{ev_1},\ldots,c^{(1)}_{ev_1+e-1}),\ldots,   
(c^{(n)}_{ev_n},\ldots,c^{(n)}_{ev_n+e-1})\!\in\!A^{e}_L$ with 
$c^{(1)}_{ev_1},\ldots,c^{(n)}_{ev_n}\!\not\in\!M_p$. 
Lipshitz's Theorem then tells us that the number of 
isolated roots $z$ of $H$ in $\C^n_p$ satisfying\\ 
$|z_1-1|_p,\ldots,|z_n-1|_p\!\leq\!\frac{1}{p}$ is no more than 
$\beta'(n,m)$. 

Since there are $(p^f-1)p^{f(e-1)}\!=\!p^d(1-\frac{1}{p^f})$ possibilities for 
each $e$-tuple $(c^{(i)}_0,\ldots,c^{(i)}_{e-1})$,  
our last observation tells us that the number of isolated roots 
$x$ of $G$ in $(L^*)^n$ satisfying 
$\ord_p x\!=\!v$ is no more than 
$p^{dn}(1-\frac{1}{p^f})^n\beta'(n,m)$. So the total number of isolated roots 
of $G$ in $(L^*)^n$ is no more than 
$\begin{pmatrix}mn\\ n+1\end{pmatrix}p^{dn}(1-\frac{1}{p^f})^n\beta'(n,m)$. 

To conclude, we simply set all possible subsets of the variables equal 
to zero (which of course never increases the number of monomial terms) 
and apply our result recursively to the resulting polynomial systems 
in fewer variables. 
We thus obtain our theorem, along with an obvious bound of 
$\gamma(n,m)\!\leq\!1+\sum^n_{i=1} \begin{pmatrix}n \\ 
i\end{pmatrix}\begin{pmatrix}mi\\ i+1\end{pmatrix}\beta'(i,m)$ for all 
$n\!\in\!\N$. \qed 

\section{Proving Our Main Global Result (Corollary \ref{cor:global})} 
\label{sec:global} 
Since $\Q$ naturally embeds in $\Q_p$ for any prime $p$, 
$K$ embeds in a degree $d$ algebraic extension, $L$, of 
$\Q_p$. So let us fix $p\!=\!2$, say. Our corollary then follows  
immediately from theorem \ref{thm:local}. 

\begin{rem}  
The following improved bound for corollary \ref{cor:global} 
follows immediately from our proof above: 
$2^{dn}(1-\frac{1}{2^{f_L}})^n\gamma(n,m)$, where $L$ is 
as in the proof. 
\end{rem} 

\section*{Acknowledgements} 
The author thanks Peter B\"urgisser for informing him of 
the important paper \cite{vandenef}, Angus McIntyre 
for suggesting this paper to Peter B\"urgisser, and 
Hendrik W.\ Lenstra, Jr.\ and Leonard Lipshitz for useful e-mail discussions. 
Special thanks 
go to an anonymous referee for graciously and elegantly pointing out some 
earlier errors of the author which were neither graceful nor elegant. 
The author is also very grateful to the architects and staff 
of {\tt MathSciNet} ({\tt http://ams.rice.edu/mathscinet/search}), 
without which, he wouldn't have found the paper \cite{lipshitz} until much 
later. 

I dedicate this paper to my god-daughter, Monica Althoff. 

\footnotesize
\bibliographystyle{acm}

\end{document}